\documentclass[draft,loadcyr]{trans2-l}

\usepackage{amsthm,amsmath,amssymb,amsfonts}

\DeclareFontFamily{OT1}{wncyr}{\hyphenchar\font45
}
\DeclareFontShape{OT1}{wncyr}{m}{n}{%
   <5> <6> <7> <8> <9> gen * wncyr
   <10> <10.95> <12> <14.4> <17.28> <20.74>  <24.88>wncyr10}{}
\DeclareFontShape{OT1}{wncyr}{m}{it}{%
   <5> <6> <7> <8> <9> gen * wncyi
   <10> <10.95> <12> <14.4> <17.28> <20.74> <24.88> wncyi10}{}
\DeclareFontShape{OT1}{wncyr}{m}{sc}{%
   <5> <6> <7> <8> <9> <10> <10.95> <12> <14.4>
   <17.28> <20.74> <24.88>wncysc10}{}
\DeclareFontShape{OT1}{wncyr}{b}{n}{%
   <5> <6> <7> <8> <9> gen * wncyb
   <10> <10.95> <12> <14.4> <17.28> <20.74> <24.88>wncyb10}{}
\input cyracc.def
\def\rus{\usefont{OT1}{wncyr}{m}{n}\cyracc\fontsize{7}{11pt}\selectfont}
\def\rusit{\usefont{OT1}{wncyr}{m}{it}\cyracc\fontsize{7}{11pt}\selectfont}

\DeclareMathSizes{9}{9}{7}{5} 

\newtheorem{theorem}{Theorem}
\newtheorem{lemma}
{Lemma}

\theoremstyle{definition}
\newtheorem{definition}
{Definition}

\theoremstyle{remark}
\newtheorem{remark}
{Remark}

\begin{document}


\newcommand{\Ce}{\mathbb C}
\newcommand{\Ha}{\mathbb H}
\newcommand{\R}{\mathbb R}
\newcommand{\OO}{\mathbb O}
\newcommand{\g}{\mathfrak g}
\newcommand{\h}{\mathfrak h}
\newcommand{\z}{\mathfrak z}
\newcommand{\ssl}{\mathfrak {sl}}
\newcommand{\ka}{\mathfrak k}
\newcommand{\el}{\mathfrak l}
\newcommand{\te}{\mathfrak t}
\newcommand{\un}{\mathfrak u}
\newcommand{\p}{\mathfrak p}
\newcommand{\al}{\mathfrak a}
\newcommand{\pt}{\!\cdot\!}
\newcommand{\pe}{\mathfrak p}
\newcommand{\gR}{{\g}^{}_{\R}}
\newcommand{\kR}{{\ka}^{}_{\R}}
\newcommand{\pR}{{\pe}^{}_{\R}}
\newcommand{\GRc}{{\rm Ad}(\g^{}_{\R})}
\newcommand{\GR}{{G}\hskip -.1mm_{\R}}
\newcommand{\cNv}{{\mathcal N}(V)}
\newcommand{\cNg}{{\mathcal N}(\g)}
\newcommand{\cNp}{{\mathcal N}(\p)}
\newcommand{\cNgr}{{\mathcal N}(\gR)}
\newcommand{\resp}{{\rm resp}.,\;\hskip -.4mm\;}
\newcommand{\blist}{\begin{enumerate}}


\title[Projective Duality
and Principal Nilpotent
Elements]{Projective Duality and Principal
Nilpotent Elements of Symmetric Pairs}

\author{Vladimir~L.~Popov}
\address{Steklov Mathematical
Institute, Russian Academy of Sciences,
Gubkina 8,  Moscow 117966, Russia}
\curraddr{} \email{popovvl@orc.ru}
\thanks{Partly
supported by Russian grant {\rus
N{SH}-123.2003.01} and by ETH Z\"urich,
Switzerland.}

\date{October 31, 2003}

\dedicatory{To A. L. Onishchik on the
occasion of his 70th birthday.\\[20pt]
\hfill{\rus To, {ch}to ne nazvano, --- ne
suwestvuet}.\\[-1.5pt]
\hfill {\rus V.~Nabokov. Prigla{sh}enie na
kazn{p1}.}}


\begin{abstract}
It is shown that projectivized irreducible
components of nilpotent cones of complex
symmetric spaces are projective self-dual
algebraic varieties. Other properties
equivalent to their projective
self-duality are found.
\end{abstract}

\maketitle

{\bf 1.}\ Let $\g$ be a semisimple
 complex Lie algebra, let $G$ be
the adjoint group of $\g$, and let $\theta
\in {\rm Aut}\,\g$ be an element of order
2. We set
\begin{equation*}
\ka:= \{x\in \g\mid \theta(x)=x\},\quad
\pe:= \{x\in \g\mid \theta(x)=-x\}.
\label{kp}
\end{equation*}
Then $\ka$ and $\pe\neq 0$, the subalgebra
$\ka$ is reductive, and $\g=\ka\oplus \pe$
is a $\mathbb Z_2$-grading of the Lie
algebra $\g$, cf.,\;e.g.,\;\cite{OV}.
Denote by $G$ the adjoint group of~$\g$.
The connected reductive algebraic subgroup
$K$ of $G$ with the Lie algebra $\ka$ is
the adjoint group of $\ka$. Denoting the
automorphism of $G$ induced by $\theta$
also by $\theta$, let $K_{\theta}$ be the
fixed point group of $\theta$. Then $K$ is
the identity component of $K_{\theta}$.

Let $\cNg$ and $\cNp$ be Zariski closed
sets of all nilpotent elements in $\g$ and
$\pe$ respecti\-ve\-ly. They are cones
(i.e.,\;stable with respect to scalar
multiplications and contain $0$). We have
\begin{equation*}
\cNp=\cNg\cap \pe. \label{intersection}
\end{equation*}
The cone $\cNg$ is irreducible, \cite{K2},
but $\cNp$, in general, is not,
cf.,\;e.g.\;\cite{Se}.

Consider the adjoint action of $G$ on
$\g$. Then  $\pe$ and  $\cNp$ are
$K_{\theta}$-stable. There are only
finitely many $G$-orbits (\resp
$K$-orbits) in $\cNg$ (\resp $\cNp$),
\cite{Dy}, \cite{K1}, \cite{KR}. Therefore
$\cNg$ (\resp every irreducible component
of $\cNp$) contains an open $G$-orbit
$\cNg_{\rm pr}$ (\resp $K$-orbit). Its
 elements are called {\it principal
nilpotent elements} of $\g$
(\resp\,$\pe$). All principal nilpotent
elements of $\pe$ constitute a single
$K_{\theta}$-orbit $\cNp_{\rm pr}$,
\cite{KR}. Hence the action of
$K_{\theta}$ on the set of irreducible
components of $\cNp$ is transitive. Remark
that there are pairs $(\g, \theta)$ for
which the intersection $\cNg_{\rm pr}\cap
\pe$ is empty.

In the sequel, if $V$ is a vector space
and $C$ is a cone in $V$, we denote by
$\mathbf P(V)$ the associated projective
space of $V$, and by $\mathbf P(C)$ the
subset of $\mathbf P(V)$ whose affine
cone~is~$C$.

The Killing form $\langle\ {,}\ \rangle$
of $\g$ (\resp its restriction $\langle\
{,}\ \rangle|{}_{\pe}$ to $\pe$) is
nondegenerate and $G$-stable (\resp
$K$-stable). Hence we may (and shall)
identify the linear spaces $\g$ and $\g^*$
(\resp $\pe$ and~$\pe^*$) by means of
$\langle\ {,}\ \rangle$ (\resp $\langle\
{,}\ \rangle|{}_{\pe}$). Then the
projective dual $\check X$ of any Zariski
closed subset $X$ of $\mathbf
P(\g)=\mathbf P(\g^*)$ (\resp $\mathbf
P(\pe)=\mathbf P(\pe^*)$),
cf.,\;e.g.,\;\cite{Ha}, becomes a Zariski
closed subset of $\mathbf P(\g)$ (\resp
$\mathbf P(\pe)$) as well. Given this, we
call $X$ {\it projective self-dual} if
$$
X=\check X.
$$

Now consider in $\mathbf P(\g)$ and
$\mathbf P(\pe)$ the Zariski closed
subsets $\mathbf P(\cNg\!)$ and $\mathbf
P(\cNp\!)$. The projective dual of
$\mathbf P(\cNg)$ was identified  in
\cite{P1}:
\begin{theorem}{\rm (\cite[Corollary 1 of
Theorem 2]{P1})}\label{null-g} The variety
$\mathbf P(\cNg\!)$ is projective
self-dual.
\end{theorem}

The goal of this note is to prove the
following
\begin{theorem}\label{null-p} Every
irreducible component of  $\mathbf
P(\cNp)$ is projective self-dual.
\end{theorem}

We also prove several other theorems
equivalent to Theorem \ref{null-p}.

Notice that the group $K_{\theta}$
transitively permutes irreducible
components of the variety $\mathbf
P(\cNp)$. Hence they are isomorphic one
another as embedded subvarieties of
$\mathbf P(\pe)$.

Some other nice geometric properties of
the variety $\mathbf P(\cNp)$ were
discovered earlier. Namely, according to
\cite{KR}, $\mathbf P(\cNp)$ is a complete
intersection in $\mathbf P(\pe)$ whose
ideal is minimally generated by $r$
homogeneous elements of the algebra
$\mathbb C[\pe]^{K}$ where $r$ is the
dimension of a Cartan subspace of $\pe$
(i.e., a maximal linear subspace of $\pe$
consisting of pairwise commuting
semisimple elements). Their degrees   are
$m_1+1,\ldots,m_r+1$, where
$m_1,\ldots,m_r$ are the exponents of the
Weyl group $W(\g,\theta)$ of the symmetric
pair $(\g,\theta)$. The dimension and the
degree of $\mathbf P(\cNp)$ are equal
respectively to $\dim \pe - r -1$ and the
order of $W(\g,\theta)$. In \cite{B},
\cite{Sl}, \cite{Se}, \cite{O}, \cite{SS},
for some symmetric pairs $(\g, \theta)$,
the generic singularities of $\mathbf
P(\cNp)$ were identified with some
simplest rational singularities. For any
$(\g, \theta)$, every irreducible
component of $\mathbf P(\cNp)$ is the
closure of some principal Hesselink
stratum of $\mathbf P(\cNp)$, \cite{He}.
Therefore Hesselink's theory yields a
desingularization of such component,
cf.\;\cite{PV}, \cite{P2} (see
also\;\cite{R}).

Theorem \ref{null-g} follows from Theorem
\ref{null-p}. Indeed, if $\g=\mathfrak
h\oplus\mathfrak h$, where $\mathfrak h$
is a semisimple complex Lie algebra, and
$\theta((y, z))=(z, y)$, then
\begin{equation}\label{reduct}
\ka=\{(y, y)\mid y\in \mathfrak h\},\quad
\p= \{(y, -y)\mid y\in \mathfrak h\}.
\end{equation}
If $H$ is the adjoint group of $\mathfrak
h$, then \eqref{reduct} implies that $K$
is isomorphic to $H$, and the $K$-module
$\p$ is isomorphic to the adjoint
$H$-module $\mathfrak h$. Identifying
these modules yields $\cNp=\mathcal
N(\mathfrak h)$.

Finally remark that in the forthcoming
paper \cite{PT} we classified (listed) all
$K$-orbits in $\mathbf P(\cNp)$ such that
their closures in $\mathbf P(\pe)$ are
projective self-dual. The $K$-orbits in
$\mathbf P(\cNp_{\rm pr})$ are not
immediately identified in this
classification. However using some extra
case by case arguments we can identify
them and thereby obtain a proof of
Theorem~\ref{null-p}. The proof given in
this note is different. It is short and
free of case by case considerations.

\bigskip

{\bf 2.} \ Our approach to proving Theorem
\ref{null-p} is based on reducing it to an
equivalent statement,
Theorem~\ref{-1-distinguished}, and then
proving the latter. This reduction is
based on the results of \cite{P1}.

To describe it, we introduce some
notation. For any subset $\mathfrak x$ of
$\g$, put
\begin{equation*}\label{pm}
{\mathfrak x}^+:=\mathfrak x\cap \ka,
\quad {\mathfrak x}^-:=\mathfrak x\cap
\pe.
\end{equation*}
If  $\mathfrak x$ and $\mathfrak y$ are
nonempty subsets of $\g$, denote by
${\mathfrak x}^\mathfrak y$ the
centralizer of $\mathfrak y$ in $\mathfrak
x$,
\begin{equation*}\label{zentralizer}
{\mathfrak x}^{\mathfrak y}:=\{x\in
\mathfrak x\mid [x, y]=0 \text{ for all }
y\in \mathfrak y\}.
\end{equation*}
If $\mathfrak x$ is a linear subspace or a
subalgebra of $\g$, then ${\mathfrak
x}^{\mathfrak y}$ has this property as
well. If $\mathfrak x$ is $\theta$-stable
and $\mathfrak y\subseteq\ka\cup\pe$, then
$ {\mathfrak x}^{\mathfrak y}$ is
$\theta$-stable.

\begin{definition} {\rm (\cite{P1})}
\label{-1} {\rm An element $x\in \cNp$ and
its $K$-orbit are called}
$(-1)$-dis\-tingui\-shed {\rm if ${\pe}^x$
contains no nonzero semisimple elements.}
\end{definition}

\begin{remark} This notion is a generalization of
the notion of distinguished nilpotent
element of a semisimple Lie algebra
introduced in the Bala--Carter theory,
\cite{BC}, cf.,\;\cite{CM}, \cite{M}.
Indeed, in the notation of \eqref{reduct},
an element of $\cNp$ is
$(-1)$-distinguished if and only if it is
distinguished as the element of $\mathcal
N(\mathfrak h)$.
\end{remark}

\begin{theorem}{\rm (\cite[Theorem 5]{P1})}
\label{duality} Let be $x$ be a nonzero
element of $\cNp$ and let
$\overline{K\cdot x}$ be the closure of
its $K$-orbit. Then the following
properties are equivalent: \blist
\item[\rm(i)] $\mathbf P(\overline{K\cdot
x})$ is projective self-dual,
\item[\rm(ii)] $x$ is
$(-1)$-distinguished.
\end{enumerate}
\end{theorem}

By Theorem \ref{duality}, Theorem
\ref{null-p} is equivalent to the
following

\begin{theorem}\label{-1-distinguished}
Every principal nilpotent element of
$\cNp$ is $(-1)$-distinguished.
\end{theorem}

\begin{remark} Since $\cNp_{\rm pr}$ is a single
$K_{\theta}$-orbit, replacing `every' with
`some' in Theorem~\ref{-1-distinguished}
 yields the equivalent
statement. The same concerns Theorems
\ref{compact}, \ref{equality} and
\ref{derived} below.
\end{remark}

Given that Theorem \ref{null-p} boils down
to Theorem ~\ref{-1-distinguished}, below
we key on proving
Theorem~\ref{-1-distinguished}.

In Subsections 6 and 7 we consider other
interesting properties of principal
nilpotent elements of $\cNp$. This yields
other statements equivalent to
Theorem~\ref{null-p}.

\bigskip

{\bf 3.}\ Since reduction of Theorem
\ref{null-p} to Theorem
\ref{-1-distinguished} is crucial for our
approach, first we sketch, for the sake of
completeness,  the proof of Theorem
\ref{duality}.

\vskip 1.3mm

\begin{proof}
The embedded tangent space to $K\cdot x$
at $x$ is $[\ka, x]$. For $X:=\mathbf
P(\overline{K\cdot x})$, the affine cone
over $\check X$ has the form
$\overline{K\cdot [\ka, x]^{\perp}}$,
where $[\ka, x]^{\perp}$ is the orthogonal
complement to $[\ka, x]$ in $\pe$ with
respect to $\langle\ {,}\
\rangle|{}_{\pe}$. By the properties of
the Killing form, $[\ka, x]^{\perp}=
{\pe}^x$.

${\rm (i)}\Rightarrow{\rm (ii)}$: Assume
that $X=\check X$. Then $ {\pe}^x\subseteq
\overline{K\cdot {\pe}^x}=\overline{K\cdot
x}\subseteq \cNp$. Thus all nonzero
elements of ${\pe}^x$ are nilpotent,
whence (ii).

${\rm (i)}\Rightarrow{\rm (ii)}$: Assume
that $x$ is $(-1)$-distinguished. Then, by
the Jordan decomposition argument,
${\pe}^x\subseteq \cNp$, whence $\check
X\subseteq \mathbf P(\cNp)$. Since there
are only finitely many $K$-orbits in
$\cNp$, the last inclusion imples that the
affine cone over $\check X$ has the form
$\overline{K\cdot y}$ for some element
$y\in \cNp$. From $x\in {\pe}^x$ and $y\in
{\pe}^y$ we deduce that $X \subseteq
\check X$ and $\check X\subseteq
\check{\check X}$. But $\check{\check
X}=X$ by the classical Biduality Theorem,
cf.\;\cite{Ha}. Whence (i).
\end{proof}

\smallskip

{\bf 4.}\ To prove Theorem
\ref{-1-distinguished} we need another
condition equivalent to
$(-1)$-dis\-tin\-guish\-ness.

In the sequel, given a real or complex
algebraic Lie algebra $\mathfrak h$, its
reductive subalgebra~$\mathfrak r$ is
called a {\it reductive Levi subalgebra}
of $\mathfrak h$ if $\mathfrak h$ is a
semidirect product of
 $\mathfrak r$ and the unipotent radical
 ${\rm rad}^{}_u\mathfrak h$ of~$\mathfrak h$.
The algebra
 $\mathfrak h/{\rm rad}^{}_u\mathfrak h$
 is called
 the {\it reductive Levi factor} of~$\mathfrak
h$.

\begin{lemma} \label{levi} Let $\mathfrak
h$ be a $\theta$-stable algebraic
subalgebra of $\g$ and let $\mathfrak r$
be a $\theta$-stable reductive Levi
subalgebra of $\mathfrak h$. Then the
following properties are equivalent:
\blist \item[\rm(i)] ${\mathfrak h}^-$
contains no nonzero semisimple elements,
\item[\rm(ii)] $\mathfrak r^-=0$.
\end{enumerate}
\end{lemma}

\begin{proof} Since $\mathfrak h$ and $\mathfrak r$
are $\theta$-stable, and $\mathfrak r$ is
a reductive Levi subalgebra of $\mathfrak
h$, we have the following direct sum
decompositions of vector spaces
\begin{equation}
\mathfrak h={\mathfrak h}^+\oplus
{\mathfrak h}^-,\quad \mathfrak
r={\mathfrak r}^+\oplus {\mathfrak
r}^-,\quad \mathfrak h=\mathfrak r\oplus
{\rm rad}_u\mathfrak h.
\label{decompositions}
\end{equation}
Since $\mathfrak h$ is $\theta$-stable,
${\rm rad}_u\mathfrak h$ is
$\theta$-stable as well. Hence we have the
decomposition
\begin{equation}
{\rm rad}_u\mathfrak h=({\rm
rad}^{}_u\mathfrak h)^+\oplus ({\rm
rad}^{}_u\mathfrak h)^-. \label{u-radical}
\end{equation}

${\rm (i)}\Rightarrow{\rm (ii)}$:  Assume
that (i) holds.
 If $\mathfrak r^-\neq
0$, then \eqref{decompositions} implies
that $\theta^{}|_{\mathfrak r}\in {\rm
Aut}\, \mathfrak r$ is an element of
order~2. Hence, by \cite{V}, there is  a
nonzero $\theta$-stable algebraic torus in
$\mathfrak r^-$. This contradicts~(i).
Whence $\mathfrak r^-=0$.

${\rm (ii)}\Rightarrow{\rm (i)}$: Assume
that (ii) holds. Then
\eqref{decompositions} implies that
\begin{equation} \mathfrak
r=\mathfrak r^+.\label{r}
\end{equation}
Plugging \eqref{r} and \eqref{u-radical}
in the last decomposition in
\eqref{decompositions}, we deduce from the
first decompo\-si\-ti\-on in
\eqref{decompositions} that
\begin{equation}
\mathfrak h^-=({\rm rad}^{}_u\mathfrak
h)^- \label{h-}.
\end{equation}
Since all elements of ${\rm
rad}^{}_u\mathfrak h$ are nilpotent,
\eqref{h-} implies (i).
\end{proof}

\medskip

Now let $e$ be a nonzero element of
$\cNp$. By Morozov's theorem, $e$ can be
embedded in an $\ssl^{}_2$-triple $\{e, h,
f\}$, i.e., an ordered triple of elements
of $\g$ satisfying the bracket relations
\begin{equation*}\label{bracket}
[h, e]=2e,\ [h, f]=-2f, \ [e, f]=h.
\end{equation*}
By \cite[Proposition 4]{KR}, we may (and
shall) assume that
\begin{equation}
h\in \ka,\quad e, f\in \pe. \label{normal}
\end{equation}

The linear span $\mathfrak s$ of $\{e, h,
f\}$ is a three dimensional subalgebra of
$\g$ isomorphic to $\ssl_2$. It is well
known (cf.,\;e.g.,\;\cite[Lemma
3.7.3]{CM}) that ${\g}^{\mathfrak s}$ is a
reductive Levi subalgebra of ${\g}^e$. By
\eqref{normal}, the subalgebras ${\g}^e$,
$\mathfrak s$ and ${\g}^{\mathfrak s}$ of
$\g$ are $\theta$-stable. Hence Lemma
\ref{levi} yields
\begin{lemma}\label{factor}
The following properties are equivalent:
\blist \item[\rm (i)] $e$ is
$(-1)$-distinguished,
\item[\rm
(ii)] ${\pe}^{\mathfrak s}=0$.
\end{enumerate}
\end{lemma}

\smallskip

{\bf 5.}\ Now we can prove Theorem {\rm
\ref{-1-distinguished}.

\begin{proof} According to the classical theory, we
may (and shall)
 fix a
$\theta$-stable real form $\gR$ of $\g$
such that
\begin{equation}
\gR=\kR\oplus\pR,\ \mbox{ where }\
\kR:={\gR}\hskip -2mm^+, \
\pR:={\gR}\hskip -2mm^-, \label{decom}
\end{equation}
is a Cartan decomposition of $\gR$,
cf.,\;e.g.,\;\cite{OV}. Let
$\al^{}_{\mathbb R}$ be a maximal abelian
subspace of $\pR$ and let $\al\subset \pe$
be its complexification. Then $\al$ is a
Cartan subspace of $\pe$,
\cite[Lemma~2]{KR}. Consider the
(restricted) root system $\Delta$ of the
pair $(\g, \al)$. Fix in $\al^{}_{\mathbb
R}$ a closed Weyl chamber $C$ of $\Delta$
and let $\Pi$ be the system of simple
roots of $\Delta$ corresponding to~$C$.

Let $e$ be a nonzero element of $\cNp$.
Include it in an $\ssl^{}_2$-triple $\{e,
h, f\}$ such that \eqref{normal} holds,
and let $\mathfrak s$ be the linear span
of $\{e, h, f\}$. By \cite[Lemma 6]{KR},
replacing $\{e, h, f\}$ with $k\cdot \{e,
h, f\}$ for an appropriate ele\-ment $k\in
K$, we may (and shall) assume that
\begin{equation}
c:=e+f\in C.\label{c}
\end{equation}
By \eqref{c}, the element $c$ lies in
  $\mathfrak s$. Therefore
\begin{equation}
{\ka}^{\mathfrak s} \subseteq
{\ka}^c,\quad {\pe}^{\mathfrak s}
\subseteq {\pe}^c. \label{incl}
\end{equation}
According to \cite[Proposition~13]{KR},
the inclusion $e\in \cNp_{\rm pr}$ is
equivalent to the property
\begin{equation}
\alpha(c)=2 \ \text{ for all } \alpha\in
\Pi.\label{condition2}
\end{equation}

Now assume that $e$ lies in $\cNp_{\rm
pr}$. Then \eqref{condition2} implies that
\begin{equation}
{\ka}^c={\ka}^{\al}, \quad
{\pe}^{c}=\al.\label{c-zentralizer}
\end{equation}

Combining \eqref{incl} and
\eqref{c-zentralizer}, we deduce from
Lemma \ref{factor} that proving the
statement of Theorem
\ref{-1-distinguished}  is equivalent to
proving the equality
\begin{equation}
{\mathfrak a}^{\mathfrak s}=0.
\label{reduction}
\end{equation}

Arguing on the contrary, assume that
\eqref{reduction} does not hold. This
means that there is a nonzero element
$z\in \mathfrak a$ commuting with every
element of $\mathfrak s$. Consider the
subalgebra $\widetilde {\mathfrak g}$ of
$\g$ generated by $\mathfrak s$ and
$\mathfrak a$. Since $\mathfrak a$ is
commutative, the element $z$ commutes with
every element of $\mathfrak a$ as well.
Therefore it commutes with every element
of the subalgebra $\widetilde {\mathfrak
g}$, i.e., belongs to its center. On the
other hand, by \cite[Proposition 23]{KR},
the algebra $\widetilde {\mathfrak g}$ is
semisimple, and hence its center is $0$.
This contradiction completes the proof of
Theorem~\ref{null-p}.
\end{proof}

\begin{remark}
It is well known that if
$x\in\cNg_{\rm pr}$, then every element of
${\g}^x$ is nilpotent, cf.\;\cite{CM}. If
$x\in\cNp_{\rm pr}$, then ${\ka}^x$, in
general, does not have this property. For
instance,
 if $\g^{}_{\mathbb R}= {\mbox{\tt
\fontsize{12pt}{0mm}\selectfont
F}}_{4(-20)}$, then the reductive Levi
factor of ${\ka}^x$ is ${\mbox{\tt
\fontsize{12pt}{0mm}\selectfont G}}_{2}$,
see \cite[Table VIII]{Do}.
\end{remark}

\smallskip

{\bf 6.}\ The decomposition \eqref{decom}
defines an $\mathbb R$-structure on the
algebraic group $G$. The identity
component of the Lie group of $\mathbb
R$-points of $G$ is the adjoint group
${\rm Ad}(\g^{}_{\R})$ of~$\gR$. We set
$$
\cNgr:=\cNg\cap \gR.
$$
\begin{definition} {\rm (\cite{PT})
An element $x\in \cNgr$ is called {\it
compact} if the reductive Levi factor of
the centralizer ${\gR}\hskip -1.5mm^x$ is
a compact Lie algebra.}
\end{definition}

Recall that there is a special bijection
between the sets of nonzero $K$-orbits in
$\cNp$ and nonzero ${\rm
Ad}(\g^{}_{\R})$-orbits in $\cNgr$,
cf.\,\cite{CM}, \cite{M}. Namely, let
$\sigma$ be the complex conjuga\-ti\-on of
$\g$ defined by $\gR$, viz.,
$$\sigma(a+ib)=a-ib, \  a, b\in \gR.$$
An $\mathfrak{sl}_2$-triple $\{e, h, f\}$
in $\g$ satisfying \eqref{normal} is
called a {\it complex Cayley triple} if
$\sigma (e)=-f$. For such a triple, set
\begin{equation*}e':=
i(-h+e+f)/2, \quad h':=e-f,\quad f':=-
i(h+e+f)/2. \label{cayley1}
\end{equation*}
Then $\{e', h', f'\}$ is an
$\mathfrak{sl}_2$-triple in $\gR$ such
that $\theta(e')=f'$. An
$\mathfrak{sl}^{}_2$-triple in $\gR$
satisfying the last property is called a
{\it real Cayley triple}. The map $\{e, h,
f\}\mapsto \{e', h', f'\}$ is a bijection
from the set of complex to the set of real
Cayley triples. The triple $\{e, h, f\}$
is called the {\it Cayley transform} of
$\{e', h', f'\}$.

Now let ${\mathcal  O}$ be a nonzero
$K$-orbit in ${\mathcal  N}(\pe)$. Then,
by \cite{KR}, there is a complex Cayley
triple $\{e, h, f\}$ in $\g$ such that
$e\in {\mathcal  O}$. Let $\{e', h', f'\}$
be the real Cayley triple in $\gR$
 such that $\{e, h, f\}$ is its Cayley
transform. Let ${\mathcal  O}'=
\GRc\!\cdot\!e'$. Then the map assigning
${\mathcal  O}'$ to $\mathcal  O$ is well
defined and establishes a bijection,
called the {\it Kostant--Sekiguchi
bijection}, between the set of nonzero
$K$-orbits in ${\mathcal  N}(\pe)$ and the
set of nonzero $\GRc$-orbits in~$\cNgr$.

\begin{theorem}{\rm
(\cite[Theorem 5]{PT})} \label{kost-sek}
Let $\mathcal  O$ be a nonzero $K$-orbit
in $\cNp$ and let $x$ be an element of the
$\GRc$-orbit in $\cNgr$ corresponding to
$\mathcal  O$ via the Kostant--Sekiguchi
bijection. Then the following pro\-perties
are equivalent: \blist \item[\rm (i)]
$\mathcal  O$ is $(-1)$-distinguished.
\item[\rm (ii)] $x$ is compact.
\end{enumerate}
\end{theorem}

Hence by Theorems \ref{kost-sek} and
\ref{-1-distinguished}, Theorem
\ref{null-p} is equivalent to the
following

\begin{theorem}\label{compact}
For every $K$-orbit $\mathcal  O$ in
$\cNp_{\rm pr}$, the elements of the
$\GRc$-orbit in $\cNgr$ corresponding to
$\mathcal  O$ via the Kostant--Sekiguchi
bijection, are compact.
\end{theorem}

\smallskip

{\bf 7.}\ Let $e$ be a nonzero element of
$ \cNp$. Fix an $\ssl_2$-triple $\{e, h,
f\}$ such that \eqref{normal} holds, and
let $\mathfrak s$ be the linear span  of
$\{e, h, f\}$. It is well known,
cf.,\;e.g.,\;\cite{CM}, \cite{M}, that
\begin{equation*}\textstyle
\g=\bigoplus_{d\in \mathbb Z} \g_d, \quad
\text{where } \ \g_d:=\{x\in \g\mid [h,
x]=dx\},\label{g}
\end{equation*}
and that
\begin{equation*}\textstyle
\mathfrak q:= \bigoplus_{d\geq
0}\g_d,\quad \mathfrak l:=\g_0,\quad
\mathfrak u= \bigoplus_{d> 0}\g_d
\end{equation*}
are respectively a parabolic subalgebra of
 $\g$ (the `Jacobson--Morozov parabolic
subalgebra' of $e$ that is actually
uniquely determined by $e$ alone), a
reductive Levi subalgebra of $\mathfrak q$
and the unipotent radical of $\mathfrak
q$. By \eqref{normal}, each $\g_d$ is
$\theta$-stable.

Since $\mathfrak s$ is $\theta$-stable,
the $\mathfrak s$-module $\g$ is  a direct
sum of $\theta$-stable simple submodules.
From this and the inclusion $e\in \g_2^-$,
we deduce by means of the known argument
based on the elementary representation
theory of $\ssl_2$
(cf.,\;e.g.,\;\cite[Lemma 8.2.1]{CM}),
that the condition $\pe^{\mathfrak s}=0$
is equivalent to the condition
\begin{equation}
\dim \g_0^-=\dim \g_2^+. \label{numeric}
\end{equation}
Using Lemma \ref{factor}, this yields the
following graded analogue of the known
criterion of distinguishness from the
Bala--Carter theory
(cf.,\;e.g.,\;\cite[Lemma 8.2.1]{CM}):

\begin{lemma}\label{02} In the notation of this subsection,
the following properties are equivalent:
\blist \item[\rm(i)] $e$ is
$(-1)$-distinguished, \item[\rm(ii)] the
equality \eqref{numeric} holds.
\end{enumerate}
\end{lemma}

Hence Theorem \ref{null-p} is equivalent
to the following
\begin{theorem}\label{equality}
The equality \eqref{numeric} holds for
every element $e\in \cNp_{\rm pr}$
satisfying~\eqref{normal}.
\end{theorem}

Now assume that the element $e$ is {\it
even}, i.e., $\g_d=0$ for every odd $d$.
Then we have the following graded analogue
of another known criterion of
distinguishness from the Bala--Carter
theory (cf.,\;e.g.,\;\cite[Lemma
8.2.6]{CM}):
\begin{lemma}\label{even} In the notation of
this subsection, if $e$ is even, then the
following properties are equivalent:
\blist \item[\rm(i)] $e$ is
$(-1)$-distinguished, \item[\rm(ii)] $\dim
{\mathfrak l}^-=\dim {\mathfrak
u}^+/[\mathfrak u, \mathfrak u]^+$.
\end{enumerate}
\end{lemma}

\begin{proof} This is deduced from Lemma \ref{02}
using the argument analogous to the one
used in the non-graded situation,
cf.,\;e.g.,\;\cite[Theorem
 8.2.6]{CM}.
 \end{proof}

\begin{remark}It is well known
(cf.,\;e.g.,\;\cite[Theorem 8.2.3]{CM})
that every distinguished element of $\cNg$
is even. In contrast to this, there are
$(-1)$-distinguished elements in $\cNp$
that are not even, \cite{PT}.
\end{remark}

Since, by \cite[Theorem 4]{KR}, the
elements of $\cNp_{\rm pr}$ are even,
Lemma \ref{even} yields that Theorem
\ref{null-p} is equivalent to the
following

\begin{theorem} \label{derived}
The equality  $\dim {\mathfrak l}^-=\dim
{\mathfrak u}^+/[\mathfrak u, \mathfrak
u]^+$ holds for every element $e\in
\cNp_{\rm pr}$ satisfying \eqref{normal}.
\end{theorem}

\smallskip

{\bf 8.}\ In \cite{N}, it was made an
attempt to develop an analogue of the
Bala--Carter theory for nilpotent orbits
in real semisimple Lie algebras. The
Kostant--Sekiguchi bijection reduces this
to finding an analogue of the Bala--Carter
theory for $K$-orbits in $\cNp$. The
theory developed in \cite{N} is based on
the notion of {\it noticed} nilpotent
element. In the notation of Subsection 4,
it is an element $e\in \cNp$ characterized
by the property $\ka^{\mathfrak s}=0$.
Recall that in contrast to this, our
$(-1)$-distinguished element is
characterized by the property
$\pe^{\mathfrak s}=0$ (see
Lemma~\ref{factor}).

The theory developed in \cite{N} does not
have some features that might be expected
from a natural analogue of the
Bala--Carter theory. As Theorem
\ref{duality} shows, the geometric
counterpart of distinguisheness of a
nonzero element $x\in \cNg$ in the
Bala--Carter theory is projective
self-duality of the variety $\mathbf
P(\overline{G\cdot x})$, and in the
contents of symmetric pairs $(\g,
\theta)$, the algebraic counterpart of
projective self-duality of $\mathbf
P(\overline{K\cdot x})$ for $x\in \cNp$ is
that $x$ is $(-1)$-distinguished, not that
$x$ is noticed.

We believe that the results of \cite{P1},
this note and \cite{PT} provide an
evidence that if a natural analogue of the
Bala--Carter theory for symmetric pairs
exists, the notion of $(-1)$-distinguished
element should play a key role in it,
analogous to that of distinguished element
in the original Bala--Carter theory.

\def\rus{\usefont{OT1}{wncyr}{m}{n}\cyracc\fontsize{8}{11pt}\selectfont}
\def\rusit{\usefont{OT1}{wncyr}{m}{it}\cyracc\fontsize{8}{11pt}\selectfont}

\bibliographystyle{amsalpha}

\end{document}